\input amstex
\input amsppt.sty

\magnification1200
\hsize13cm
\vsize19cm

\TagsOnRight
\NoBlackBoxes

\def\RoZeAA{7}
\def\MacdAC{6}
\def\LindAA{5}
\def\LaWeAB{4}
\def\GeViAB{3}
\def\GeViAA{2}
\def\FultAC{1}

\catcode`\@=11
\def\dddot#1{\vbox{\ialign{##\crcr
      .\hskip-.5pt.\hskip-.5pt.\crcr\noalign{\kern1.5\p@\nointerlineskip}
      $\hfil\displaystyle{#1}\hfil$\crcr}}}

\newif\iftab@\tab@false
\newif\ifvtab@\vtab@false
\def\tab{\bgroup\tab@true\vtab@false\vst@bfalse\Strich@false%
   \def\\{\global\hline@@false%
     \ifhline@\global\hline@false\global\hline@@true\fi\cr}
   \edef\l@{\the\leftskip}\ialign\bgroup\hskip\l@##\hfil&&##\hfil\cr}
\def\endtab{\cr\egroup\egroup}
\def\vtab{\vtop\bgroup\vst@bfalse\vtab@true\tab@true\Strich@false%
   \bgroup\def\\{\cr}\ialign\bgroup&##\hfil\cr}
\def\endvtab{\cr\egroup\egroup\egroup}
\def\stab{\D@cke0.5pt\null 
 \bgroup\tab@true\vtab@false\vst@bfalse\Strich@true\Let@@\vspace@
 \normalbaselines\offinterlineskip
  \openup\spreadmlines@
 \edef\l@{\the\leftskip}\ialign
 \bgroup\hskip\l@##\hfil&&##\hfil\crcr}
\def\endstab{\crcr\egroup
 \egroup}
\newif\ifvst@b\vst@bfalse
\def\vstab{\D@cke0.5pt\null
 \vtop\bgroup\tab@true\vtab@false\vst@btrue\Strich@true\bgroup\Let@@\vspace@
 \normalbaselines\offinterlineskip
  \openup\spreadmlines@\bgroup}
\def\endvstab{\crcr\egroup\egroup
 \egroup\tab@false\Strich@false}

\newdimen\htstrut@
\htstrut@8.5\p@
\newdimen\htStrut@
\htStrut@12\p@
\newdimen\dpstrut@
\dpstrut@3.5\p@
\newdimen\dpStrut@
\dpStrut@3.5\p@
\def\openup{\afterassignment\@penup\dimen@=}
\def\@penup{\advance\lineskip\dimen@
  \advance\baselineskip\dimen@
  \advance\lineskiplimit\dimen@
  \divide\dimen@ by2
  \advance\htstrut@\dimen@
  \advance\htStrut@\dimen@
  \advance\dpstrut@\dimen@
  \advance\dpStrut@\dimen@}
\def\Let@@{\relax%
    \def\\{\global\hline@@false%
     \ifhline@\global\hline@false\global\hline@@true\fi\cr}%
    \iffalse}\fi}
\def\matrix{\null\,\vcenter\bgroup
 \tab@false\vtab@false\vst@bfalse\Strich@false\Let@@\vspace@
 \normalbaselines\openup\spreadmlines@\ialign
 \bgroup\hfil$\m@th##$\hfil&&\quad\hfil$\m@th##$\hfil\crcr
 \Mathstrut@\crcr\noalign{\kern-\baselineskip}}
\def\endmatrix{\crcr\Mathstrut@\crcr\noalign{\kern-\baselineskip}\egroup
 \egroup\,}
\def\smatrix{\D@cke0.5pt\null\,
 \vcenter\bgroup\tab@false\vtab@false\vst@bfalse\Strich@true\Let@@\vspace@
 \normalbaselines\offinterlineskip
  \openup\spreadmlines@\ialign
 \bgroup\hfil$\m@th##$\hfil&&\quad\hfil$\m@th##$\hfil\crcr}
\def\endsmatrix{\crcr\egroup
 \egroup\,\Strich@false}
\newdimen\D@cke
\def\Dicke#1{\global\D@cke#1}
\newtoks\tabs@\tabs@{&}
\newif\ifStrich@\Strich@false
\newif\iff@rst

\def\Stricherr@{\iftab@\ifvtab@\errmessage{\noexpand\s not allowed
     here. Use \noexpand\vstab!}%
  \else\errmessage{\noexpand\s not allowed here. Use \noexpand\stab!}%
  \fi\else\errmessage{\noexpand\s not allowed
     here. Use \noexpand\smatrix!}\fi}
\def\format{\ifvst@b\else\crcr\fi\egroup\iffalse{\fi\ifnum`}=0 \fi\format@}
\def\format@#1\\{\def\preamble@{#1}%
 \def\Str@chfehlt##1{\ifx##1\s\Stricherr@\fi\ifx##1\\\let\Next\relax%
   \else\let\Next\Str@chfehlt\fi\Next}%
 \def\c{\hfil\noexpand\ifhline@@\hbox{\vrule height\htStrut@%
   depth\dpstrut@ width\z@}\noexpand\fi%
   \ifStrich@\hbox{\vrule height\htstrut@ depth\dpstrut@ width\z@}%
   \fi\iftab@\else$\m@th\fi\the\hashtoks@\iftab@\else$\fi\hfil}%
 \def\r{\hfil\noexpand\ifhline@@\hbox{\vrule height\htStrut@%
   depth\dpstrut@ width\z@}\noexpand\fi%
   \ifStrich@\hbox{\vrule height\htstrut@ depth\dpstrut@ width\z@}%
   \fi\iftab@\else$\m@th\fi\the\hashtoks@\iftab@\else$\fi}%
 \def\l{\noexpand\ifhline@@\hbox{\vrule height\htStrut@%
   depth\dpstrut@ width\z@}\noexpand\fi%
   \ifStrich@\hbox{\vrule height\htstrut@ depth\dpstrut@ width\z@}%
   \fi\iftab@\else$\m@th\fi\the\hashtoks@\iftab@\else$\fi\hfil}%
 \def\s{\ifStrich@\ \the\tabs@\vrule width\D@cke\the\hashtoks@%
          \fi\the\tabs@\ }%
 \def\sa{\ifStrich@\vrule width\D@cke\the\hashtoks@%
            \the\tabs@\ %
            \fi}%
 \def\se{\ifStrich@\ \the\tabs@\vrule width\D@cke\the\hashtoks@\fi}%
 \def\cd{\hfil\noexpand\ifhline@@\hbox{\vrule height\htStrut@%
   depth\dpstrut@ width\z@}\noexpand\fi%
   \ifStrich@\hbox{\vrule height\htstrut@ depth\dpstrut@ width\z@}%
   \fi$\dsize\m@th\the\hashtoks@$\hfil}%
 \def\rd{\hfil\noexpand\ifhline@@\hbox{\vrule height\htStrut@%
   depth\dpstrut@ width\z@}\noexpand\fi%
   \ifStrich@\hbox{\vrule height\htstrut@ depth\dpstrut@ width\z@}%
   \fi$\dsize\m@th\the\hashtoks@$}%
 \def\ld{\noexpand\ifhline@@\hbox{\vrule height\htStrut@%
   depth\dpstrut@ width\z@}\noexpand\fi%
   \ifStrich@\hbox{\vrule height\htstrut@ depth\dpstrut@ width\z@}%
   \fi$\dsize\m@th\the\hashtoks@$\hfil}%
 \ifStrich@\else\Str@chfehlt#1\\\fi%
 \setbox\z@\hbox{\xdef\Preamble@{\preamble@}}\ifnum`{=0 \fi\iffalse}\fi
 \ialign\bgroup\span\Preamble@\crcr}
\newif\ifhline@\hline@false
\newif\ifhline@@\hline@@false
\def\hlinefor#1{\multispan@{\strip@#1 }\leaders\hrule height\D@cke\hfill%
    \global\hline@true\ignorespaces}
\catcode`\@=13

\catcode`\@=11
\font\tenln    = line10
\font\tenlnw   = linew10

\newskip\Einheit \Einheit=0.5cm
\newcount\xcoord \newcount\ycoord
\newdimen\xdim \newdimen\ydim \newdimen\PfadD@cke \newdimen\Pfadd@cke

\newcount\@tempcnta
\newcount\@tempcntb

\newdimen\@tempdima
\newdimen\@tempdimb

\newdimen\@wholewidth
\newdimen\@halfwidth

\newcount\@xarg
\newcount\@yarg
\newcount\@yyarg
\newbox\@linechar
\newbox\@tempboxa
\newdimen\@linelen
\newdimen\@clnwd
\newdimen\@clnht

\newif\if@negarg

\def\@whilenoop#1{}
\def\@whiledim#1\do #2{\ifdim #1\relax#2\@iwhiledim{#1\relax#2}\fi}
\def\@iwhiledim#1{\ifdim #1\let\@nextwhile=\@iwhiledim
        \else\let\@nextwhile=\@whilenoop\fi\@nextwhile{#1}}

\def\@whileswnoop#1\fi{}
\def\@whilesw#1\fi#2{#1#2\@iwhilesw{#1#2}\fi\fi}
\def\@iwhilesw#1\fi{#1\let\@nextwhile=\@iwhilesw
         \else\let\@nextwhile=\@whileswnoop\fi\@nextwhile{#1}\fi}

\def\thinlines{\let\@linefnt\tenln \let\@circlefnt\tencirc
  \@wholewidth\fontdimen8\tenln \@halfwidth .5\@wholewidth}
\def\thicklines{\let\@linefnt\tenlnw \let\@circlefnt\tencircw
  \@wholewidth\fontdimen8\tenlnw \@halfwidth .5\@wholewidth}
\thinlines

\PfadD@cke1pt \Pfadd@cke0.5pt
\def\PfadDicke#1{\PfadD@cke#1 \divide\PfadD@cke by2 \Pfadd@cke\PfadD@cke \multiply\PfadD@cke by2}
\long\def\LOOP#1\REPEAT{\def\BODY{#1}\ITERATE}
\def\ITERATE{\BODY \let\next\ITERATE \else\let\next\relax\fi \next}
\let\REPEAT=\fi
\def\Punkt{\hbox{\raise-2pt\hbox to0pt{\hss$\ssize\bullet$\hss}}}
\def\DuennPunkt(#1,#2){\unskip
  \raise#2 \Einheit\hbox to0pt{\hskip#1 \Einheit
          \raise-2.5pt\hbox to0pt{\hss$\bullet$\hss}\hss}}
\def\NormalPunkt(#1,#2){\unskip
  \raise#2 \Einheit\hbox to0pt{\hskip#1 \Einheit
          \raise-3pt\hbox to0pt{\hss\twelvepoint$\bullet$\hss}\hss}}
\def\DickPunkt(#1,#2){\unskip
  \raise#2 \Einheit\hbox to0pt{\hskip#1 \Einheit
          \raise-4pt\hbox to0pt{\hss\fourteenpoint$\bullet$\hss}\hss}}
\def\Kreis(#1,#2){\unskip
  \raise#2 \Einheit\hbox to0pt{\hskip#1 \Einheit
          \raise-4pt\hbox to0pt{\hss\fourteenpoint$\circ$\hss}\hss}}

\def\Line@(#1,#2)#3{\@xarg #1\relax \@yarg #2\relax
\@linelen=#3\Einheit
\ifnum\@xarg =0 \@vline
  \else \ifnum\@yarg =0 \@hline \else \@sline\fi
\fi}

\def\@sline{\ifnum\@xarg< 0 \@negargtrue \@xarg -\@xarg \@yyarg -\@yarg
  \else \@negargfalse \@yyarg \@yarg \fi
\ifnum \@yyarg >0 \@tempcnta\@yyarg \else \@tempcnta -\@yyarg \fi
\ifnum\@tempcnta>6 \@badlinearg\@tempcnta0 \fi
\ifnum\@xarg>6 \@badlinearg\@xarg 1 \fi
\setbox\@linechar\hbox{\@linefnt\@getlinechar(\@xarg,\@yyarg)}%
\ifnum \@yarg >0 \let\@upordown\raise \@clnht\z@
   \else\let\@upordown\lower \@clnht \ht\@linechar\fi
\@clnwd=\wd\@linechar
\if@negarg \hskip -\wd\@linechar \def\@tempa{\hskip -2\wd\@linechar}\else
     \let\@tempa\relax \fi
\@whiledim \@clnwd <\@linelen \do
  {\@upordown\@clnht\copy\@linechar
   \@tempa
   \advance\@clnht \ht\@linechar
   \advance\@clnwd \wd\@linechar}%
\advance\@clnht -\ht\@linechar
\advance\@clnwd -\wd\@linechar
\@tempdima\@linelen\advance\@tempdima -\@clnwd
\@tempdimb\@tempdima\advance\@tempdimb -\wd\@linechar
\if@negarg \hskip -\@tempdimb \else \hskip \@tempdimb \fi
\multiply\@tempdima \@m
\@tempcnta \@tempdima \@tempdima \wd\@linechar \divide\@tempcnta \@tempdima
\@tempdima \ht\@linechar \multiply\@tempdima \@tempcnta
\divide\@tempdima \@m
\advance\@clnht \@tempdima
\ifdim \@linelen <\wd\@linechar
   \hskip \wd\@linechar
  \else\@upordown\@clnht\copy\@linechar\fi}

\def\@hline{\ifnum \@xarg <0 \hskip -\@linelen \fi
\vrule height\Pfadd@cke width \@linelen depth\Pfadd@cke
\ifnum \@xarg <0 \hskip -\@linelen \fi}

\def\@getlinechar(#1,#2){\@tempcnta#1\relax\multiply\@tempcnta 8
\advance\@tempcnta -9 \ifnum #2>0 \advance\@tempcnta #2\relax\else
\advance\@tempcnta -#2\relax\advance\@tempcnta 64 \fi
\char\@tempcnta}

\def\Vektor(#1,#2)#3(#4,#5){\unskip\leavevmode
  \xcoord#4\relax \ycoord#5\relax
      \raise\ycoord \Einheit\hbox to0pt{\hskip\xcoord \Einheit
         \Vector@(#1,#2){#3}\hss}}

\def\Vector@(#1,#2)#3{\@xarg #1\relax \@yarg #2\relax
\@tempcnta \ifnum\@xarg<0 -\@xarg\else\@xarg\fi
\ifnum\@tempcnta<5\relax
\@linelen=#3\Einheit
\ifnum\@xarg =0 \@vvector
  \else \ifnum\@yarg =0 \@hvector \else \@svector\fi
\fi
\else\@badlinearg\fi}

\def\@hvector{\@hline\hbox to 0pt{\@linefnt
\ifnum \@xarg <0 \@getlarrow(1,0)\hss\else
    \hss\@getrarrow(1,0)\fi}}

\def\@vvector{\ifnum \@yarg <0 \@downvector \else \@upvector \fi}

\def\@svector{\@sline
\@tempcnta\@yarg \ifnum\@tempcnta <0 \@tempcnta=-\@tempcnta\fi
\ifnum\@tempcnta <5
  \hskip -\wd\@linechar
  \@upordown\@clnht \hbox{\@linefnt  \if@negarg
  \@getlarrow(\@xarg,\@yyarg) \else \@getrarrow(\@xarg,\@yyarg) \fi}%
\else\@badlinearg\fi}

\def\@upline{\hbox to \z@{\hskip -.5\Pfadd@cke \vrule width \Pfadd@cke
   height \@linelen depth \z@\hss}}

\def\@downline{\hbox to \z@{\hskip -.5\Pfadd@cke \vrule width \Pfadd@cke
   height \z@ depth \@linelen \hss}}

\def\@upvector{\@upline\setbox\@tempboxa\hbox{\@linefnt\char'66}\raise
     \@linelen \hbox to\z@{\lower \ht\@tempboxa\box\@tempboxa\hss}}

\def\@downvector{\@downline\lower \@linelen
      \hbox to \z@{\@linefnt\char'77\hss}}

\def\@getlarrow(#1,#2){\ifnum #2 =\z@ \@tempcnta='33\else
\@tempcnta=#1\relax\multiply\@tempcnta \sixt@@n \advance\@tempcnta
-9 \@tempcntb=#2\relax\multiply\@tempcntb \tw@
\ifnum \@tempcntb >0 \advance\@tempcnta \@tempcntb\relax
\else\advance\@tempcnta -\@tempcntb\advance\@tempcnta 64
\fi\fi\char\@tempcnta}

\def\@getrarrow(#1,#2){\@tempcntb=#2\relax
\ifnum\@tempcntb < 0 \@tempcntb=-\@tempcntb\relax\fi
\ifcase \@tempcntb\relax \@tempcnta='55 \or
\ifnum #1<3 \@tempcnta=#1\relax\multiply\@tempcnta
24 \advance\@tempcnta -6 \else \ifnum #1=3 \@tempcnta=49
\else\@tempcnta=58 \fi\fi\or
\ifnum #1<3 \@tempcnta=#1\relax\multiply\@tempcnta
24 \advance\@tempcnta -3 \else \@tempcnta=51\fi\or
\@tempcnta=#1\relax\multiply\@tempcnta
\sixt@@n \advance\@tempcnta -\tw@ \else
\@tempcnta=#1\relax\multiply\@tempcnta
\sixt@@n \advance\@tempcnta 7 \fi\ifnum #2<0 \advance\@tempcnta 64 \fi
\char\@tempcnta}

\def\Diagonale(#1,#2)#3{\unskip\leavevmode
  \xcoord#1\relax \ycoord#2\relax
      \raise\ycoord \Einheit\hbox to0pt{\hskip\xcoord \Einheit
         \Line@(1,1){#3}\hss}}
\def\AntiDiagonale(#1,#2)#3{\unskip\leavevmode
  \xcoord#1\relax \ycoord#2\relax 
      \raise\ycoord \Einheit\hbox to0pt{\hskip\xcoord \Einheit
         \Line@(1,-1){#3}\hss}}
\def\Pfad(#1,#2),#3\endPfad{\unskip\leavevmode
  \xcoord#1 \ycoord#2 \thicklines\ZeichnePfad#3\endPfad\thinlines}
\def\ZeichnePfad#1{\ifx#1\endPfad\let\next\relax
  \else\let\next\ZeichnePfad
    \ifnum#1=1
      \raise\ycoord \Einheit\hbox to0pt{\hskip\xcoord \Einheit
         \vrule height\Pfadd@cke width1 \Einheit depth\Pfadd@cke\hss}%
      \advance\xcoord by 1
    \else\ifnum#1=2
      \raise\ycoord \Einheit\hbox to0pt{\hskip\xcoord \Einheit
        \hbox{\hskip-\PfadD@cke\vrule height1 \Einheit width\PfadD@cke depth0pt}\hss}%
      \advance\ycoord by 1
    \else\ifnum#1=3
      \raise\ycoord \Einheit\hbox to0pt{\hskip\xcoord \Einheit
         \Line@(1,1){1}\hss}
      \advance\xcoord by 1
      \advance\ycoord by 1
    \else\ifnum#1=4
      \raise\ycoord \Einheit\hbox to0pt{\hskip\xcoord \Einheit
         \Line@(1,-1){1}\hss}
      \advance\xcoord by 1
      \advance\ycoord by -1
    \else\ifnum#1=5
      \advance\xcoord by -1
      \raise\ycoord \Einheit\hbox to0pt{\hskip\xcoord \Einheit
         \vrule height\Pfadd@cke width1 \Einheit depth\Pfadd@cke\hss}%
    \else\ifnum#1=6
      \advance\ycoord by -1
      \raise\ycoord \Einheit\hbox to0pt{\hskip\xcoord \Einheit
        \hbox{\hskip-\PfadD@cke\vrule height1 \Einheit width\PfadD@cke depth0pt}\hss}%
    \else\ifnum#1=7
      \advance\xcoord by -1
      \advance\ycoord by -1
      \raise\ycoord \Einheit\hbox to0pt{\hskip\xcoord \Einheit
         \Line@(1,1){1}\hss}
    \else\ifnum#1=8
      \advance\xcoord by -1
      \advance\ycoord by +1
      \raise\ycoord \Einheit\hbox to0pt{\hskip\xcoord \Einheit
         \Line@(1,-1){1}\hss}
    \fi\fi\fi\fi
    \fi\fi\fi\fi
  \fi\next}
\def\hSSchritt{\leavevmode\raise-.4pt\hbox to0pt{\hss.\hss}\hskip.2\Einheit
  \raise-.4pt\hbox to0pt{\hss.\hss}\hskip.2\Einheit
  \raise-.4pt\hbox to0pt{\hss.\hss}\hskip.2\Einheit
  \raise-.4pt\hbox to0pt{\hss.\hss}\hskip.2\Einheit
  \raise-.4pt\hbox to0pt{\hss.\hss}\hskip.2\Einheit}
\def\vSSchritt{\vbox{\baselineskip.2\Einheit\lineskiplimit0pt
\hbox{.}\hbox{.}\hbox{.}\hbox{.}\hbox{.}}}
\def\DSSchritt{\leavevmode\raise-.4pt\hbox to0pt{%
  \hbox to0pt{\hss.\hss}\hskip.2\Einheit
  \raise.2\Einheit\hbox to0pt{\hss.\hss}\hskip.2\Einheit
  \raise.4\Einheit\hbox to0pt{\hss.\hss}\hskip.2\Einheit
  \raise.6\Einheit\hbox to0pt{\hss.\hss}\hskip.2\Einheit
  \raise.8\Einheit\hbox to0pt{\hss.\hss}\hss}}
\def\dSSchritt{\leavevmode\raise-.4pt\hbox to0pt{%
  \hbox to0pt{\hss.\hss}\hskip.2\Einheit
  \raise-.2\Einheit\hbox to0pt{\hss.\hss}\hskip.2\Einheit
  \raise-.4\Einheit\hbox to0pt{\hss.\hss}\hskip.2\Einheit
  \raise-.6\Einheit\hbox to0pt{\hss.\hss}\hskip.2\Einheit
  \raise-.8\Einheit\hbox to0pt{\hss.\hss}\hss}}
\def\SPfad(#1,#2),#3\endSPfad{\unskip\leavevmode
  \xcoord#1 \ycoord#2 \ZeichneSPfad#3\endSPfad}
\def\ZeichneSPfad#1{\ifx#1\endSPfad\let\next\relax
  \else\let\next\ZeichneSPfad
    \ifnum#1=1
      \raise\ycoord \Einheit\hbox to0pt{\hskip\xcoord \Einheit
         \hSSchritt\hss}%
      \advance\xcoord by 1
    \else\ifnum#1=2
      \raise\ycoord \Einheit\hbox to0pt{\hskip\xcoord \Einheit
        \hbox{\hskip-2pt \vSSchritt}\hss}%
      \advance\ycoord by 1
    \else\ifnum#1=3
      \raise\ycoord \Einheit\hbox to0pt{\hskip\xcoord \Einheit
         \DSSchritt\hss}
      \advance\xcoord by 1
      \advance\ycoord by 1
    \else\ifnum#1=4
      \raise\ycoord \Einheit\hbox to0pt{\hskip\xcoord \Einheit
         \dSSchritt\hss}
      \advance\xcoord by 1
      \advance\ycoord by -1
    \else\ifnum#1=5
      \advance\xcoord by -1
      \raise\ycoord \Einheit\hbox to0pt{\hskip\xcoord \Einheit
         \hSSchritt\hss}%
    \else\ifnum#1=6
      \advance\ycoord by -1
      \raise\ycoord \Einheit\hbox to0pt{\hskip\xcoord \Einheit
        \hbox{\hskip-2pt \vSSchritt}\hss}%
    \else\ifnum#1=7
      \advance\xcoord by -1
      \advance\ycoord by -1
      \raise\ycoord \Einheit\hbox to0pt{\hskip\xcoord \Einheit
         \DSSchritt\hss}
    \else\ifnum#1=8
      \advance\xcoord by -1
      \advance\ycoord by 1
      \raise\ycoord \Einheit\hbox to0pt{\hskip\xcoord \Einheit
         \dSSchritt\hss}
    \fi\fi\fi\fi
    \fi\fi\fi\fi
  \fi\next}
\def\Koordinatenachsen(#1,#2){\unskip
 \hbox to0pt{\hskip-.5pt\vrule height#2 \Einheit width.5pt depth1 \Einheit}%
 \hbox to0pt{\hskip-1 \Einheit \xcoord#1 \advance\xcoord by1
    \vrule height0.25pt width\xcoord \Einheit depth0.25pt\hss}}
\def\Koordinatenachsen(#1,#2)(#3,#4){\unskip
 \hbox to0pt{\hskip-.5pt \ycoord-#4 \advance\ycoord by1
    \vrule height#2 \Einheit width.5pt depth\ycoord \Einheit}%
 \hbox to0pt{\hskip-1 \Einheit \hskip#3\Einheit 
    \xcoord#1 \advance\xcoord by1 \advance\xcoord by-#3 
    \vrule height0.25pt width\xcoord \Einheit depth0.25pt\hss}}
\def\Gitter(#1,#2){\unskip \xcoord0 \ycoord0 \leavevmode
  \LOOP\ifnum\ycoord<#2
    \loop\ifnum\xcoord<#1
      \raise\ycoord \Einheit\hbox to0pt{\hskip\xcoord \Einheit\Punkt\hss}%
      \advance\xcoord by1
    \repeat
    \xcoord0
    \advance\ycoord by1
  \REPEAT}
\def\Gitter(#1,#2)(#3,#4){\unskip \xcoord#3 \ycoord#4 \leavevmode
  \LOOP\ifnum\ycoord<#2
    \loop\ifnum\xcoord<#1
      \raise\ycoord \Einheit\hbox to0pt{\hskip\xcoord \Einheit\Punkt\hss}%
      \advance\xcoord by1
    \repeat
    \xcoord#3
    \advance\ycoord by1
  \REPEAT}
\def\Label#1#2(#3,#4){\unskip \xdim#3 \Einheit \ydim#4 \Einheit
  \def\lo{\advance\xdim by-.5 \Einheit \advance\ydim by.5 \Einheit}%
  \def\llo{\advance\xdim by-.25cm \advance\ydim by.5 \Einheit}%
  \def\loo{\advance\xdim by-.5 \Einheit \advance\ydim by.25cm}%
  \def\o{\advance\ydim by.25cm}%
  \def\ro{\advance\xdim by.5 \Einheit \advance\ydim by.5 \Einheit}%
  \def\rro{\advance\xdim by.25cm \advance\ydim by.5 \Einheit}%
  \def\roo{\advance\xdim by.5 \Einheit \advance\ydim by.25cm}%
  \def\l{\advance\xdim by-.30cm}%
  \def\r{\advance\xdim by.30cm}%
  \def\lu{\advance\xdim by-.5 \Einheit \advance\ydim by-.6 \Einheit}%
  \def\llu{\advance\xdim by-.25cm \advance\ydim by-.6 \Einheit}%
  \def\luu{\advance\xdim by-.5 \Einheit \advance\ydim by-.30cm}%
  \def\u{\advance\ydim by-.30cm}%
  \def\ru{\advance\xdim by.5 \Einheit \advance\ydim by-.6 \Einheit}%
  \def\rru{\advance\xdim by.25cm \advance\ydim by-.6 \Einheit}%
  \def\ruu{\advance\xdim by.5 \Einheit \advance\ydim by-.30cm}%
  #1\raise\ydim\hbox to0pt{\hskip\xdim
     \vbox to0pt{\vss\hbox to0pt{\hss$#2$\hss}\vss}\hss}%
}
\catcode`\@=13

\catcode`\@=11
\font@\twelverm=cmr10 scaled\magstep1
\font@\twelveit=cmti10 scaled\magstep1
\font@\twelvebf=cmbx10 scaled\magstep1
\font@\twelvei=cmmi10 scaled\magstep1
\font@\twelvesy=cmsy10 scaled\magstep1
\font@\twelveex=cmex10 scaled\magstep1

\newtoks\twelvepoint@
\def\twelvepoint{\normalbaselineskip15\p@
 \abovedisplayskip15\p@ plus3.6\p@ minus10.8\p@
 \belowdisplayskip\abovedisplayskip
 \abovedisplayshortskip\z@ plus3.6\p@
 \belowdisplayshortskip8.4\p@ plus3.6\p@ minus4.8\p@
 \textonlyfont@\rm\twelverm \textonlyfont@\it\twelveit
 \textonlyfont@\sl\twelvesl \textonlyfont@\bf\twelvebf
 \textonlyfont@\smc\twelvesmc \textonlyfont@\tt\twelvett
%
 \ifsyntax@ \def\big##1{{\hbox{$\left##1\right.$}}}%
  \let\Big\big \let\bigg\big \let\Bigg\big
 \else
  \textfont\z@=\twelverm  \scriptfont\z@=\tenrm  \scriptscriptfont\z@=\sevenrm
  \textfont\@ne=\twelvei  \scriptfont\@ne=\teni  \scriptscriptfont\@ne=\seveni
  \textfont\tw@=\twelvesy \scriptfont\tw@=\tensy \scriptscriptfont\tw@=\sevensy
  \textfont\thr@@=\twelveex \scriptfont\thr@@=\tenex
        \scriptscriptfont\thr@@=\tenex
  \textfont\itfam=\twelveit \scriptfont\itfam=\tenit
        \scriptscriptfont\itfam=\tenit
  \textfont\bffam=\twelvebf \scriptfont\bffam=\tenbf
        \scriptscriptfont\bffam=\sevenbf
  \setbox\strutbox\hbox{\vrule height10.2\p@ depth4.2\p@ width\z@}%
  \setbox\strutbox@\hbox{\lower.6\normallineskiplimit\vbox{%
        \kern-\normallineskiplimit\copy\strutbox}}%
 \setbox\z@\vbox{\hbox{$($}\kern\z@}\bigsize@=1.4\ht\z@
 \fi
 \normalbaselines\rm\ex@.2326ex\jot3.6\ex@\the\twelvepoint@}

\font@\fourteenrm=cmr10 scaled\magstep2
\font@\fourteenit=cmti10 scaled\magstep2
\font@\fourteensl=cmsl10 scaled\magstep2
\font@\fourteensmc=cmcsc10 scaled\magstep2
\font@\fourteentt=cmtt10 scaled\magstep2
\font@\fourteenbf=cmbx10 scaled\magstep2
\font@\fourteeni=cmmi10 scaled\magstep2
\font@\fourteensy=cmsy10 scaled\magstep2
\font@\fourteenex=cmex10 scaled\magstep2
\font@\fourteenmsa=msam10 scaled\magstep2
\font@\fourteeneufm=eufm10 scaled\magstep2
\font@\fourteenmsb=msbm10 scaled\magstep2
\newtoks\fourteenpoint@
\def\fourteenpoint{\normalbaselineskip15\p@
 \abovedisplayskip18\p@ plus4.3\p@ minus12.9\p@
 \belowdisplayskip\abovedisplayskip
 \abovedisplayshortskip\z@ plus4.3\p@
 \belowdisplayshortskip10.1\p@ plus4.3\p@ minus5.8\p@
 \textonlyfont@\rm\fourteenrm \textonlyfont@\it\fourteenit
 \textonlyfont@\sl\fourteensl \textonlyfont@\bf\fourteenbf
 \textonlyfont@\smc\fourteensmc \textonlyfont@\tt\fourteentt
%
 \ifsyntax@ \def\big##1{{\hbox{$\left##1\right.$}}}%
  \let\Big\big \let\bigg\big \let\Bigg\big
 \else
  \textfont\z@=\fourteenrm  \scriptfont\z@=\twelverm  \scriptscriptfont\z@=\tenrm
  \textfont\@ne=\fourteeni  \scriptfont\@ne=\twelvei  \scriptscriptfont\@ne=\teni
  \textfont\tw@=\fourteensy \scriptfont\tw@=\twelvesy \scriptscriptfont\tw@=\tensy
  \textfont\thr@@=\fourteenex \scriptfont\thr@@=\twelveex
        \scriptscriptfont\thr@@=\twelveex
  \textfont\itfam=\fourteenit \scriptfont\itfam=\twelveit
        \scriptscriptfont\itfam=\twelveit
  \textfont\bffam=\fourteenbf \scriptfont\bffam=\twelvebf
        \scriptscriptfont\bffam=\tenbf
  \setbox\strutbox\hbox{\vrule height12.2\p@ depth5\p@ width\z@}%
  \setbox\strutbox@\hbox{\lower.72\normallineskiplimit\vbox{%
        \kern-\normallineskiplimit\copy\strutbox}}%
 \setbox\z@\vbox{\hbox{$($}\kern\z@}\bigsize@=1.7\ht\z@
 \fi
 \normalbaselines\rm\ex@.2326ex\jot4.3\ex@\the\fourteenpoint@}

\font@\seventeenrm=cmr10 scaled\magstep3
\font@\seventeenit=cmti10 scaled\magstep3
\font@\seventeensl=cmsl10 scaled\magstep3
\font@\seventeensmc=cmcsc10 scaled\magstep3
\font@\seventeentt=cmtt10 scaled\magstep3
\font@\seventeenbf=cmbx10 scaled\magstep3
\font@\seventeeni=cmmi10 scaled\magstep3
\font@\seventeensy=cmsy10 scaled\magstep3
\font@\seventeenex=cmex10 scaled\magstep3
\font@\seventeenmsa=msam10 scaled\magstep3
\font@\seventeeneufm=eufm10 scaled\magstep3
\font@\seventeenmsb=msbm10 scaled\magstep3
\newtoks\seventeenpoint@
\def\seventeenpoint{\normalbaselineskip18\p@
 \abovedisplayskip21.6\p@ plus5.2\p@ minus15.4\p@
 \belowdisplayskip\abovedisplayskip
 \abovedisplayshortskip\z@ plus5.2\p@
 \belowdisplayshortskip12.1\p@ plus5.2\p@ minus7\p@
 \textonlyfont@\rm\seventeenrm \textonlyfont@\it\seventeenit
 \textonlyfont@\sl\seventeensl \textonlyfont@\bf\seventeenbf
 \textonlyfont@\smc\seventeensmc \textonlyfont@\tt\seventeentt
%
 \ifsyntax@ \def\big##1{{\hbox{$\left##1\right.$}}}%
  \let\Big\big \let\bigg\big \let\Bigg\big
 \else
  \textfont\z@=\seventeenrm  \scriptfont\z@=\fourteenrm  \scriptscriptfont\z@=\twelverm
  \textfont\@ne=\seventeeni  \scriptfont\@ne=\fourteeni  \scriptscriptfont\@ne=\twelvei
  \textfont\tw@=\seventeensy \scriptfont\tw@=\fourteensy \scriptscriptfont\tw@=\twelvesy
  \textfont\thr@@=\seventeenex \scriptfont\thr@@=\fourteenex
        \scriptscriptfont\thr@@=\fourteenex
  \textfont\itfam=\seventeenit \scriptfont\itfam=\fourteenit
        \scriptscriptfont\itfam=\fourteenit
  \textfont\bffam=\seventeenbf \scriptfont\bffam=\fourteenbf
        \scriptscriptfont\bffam=\twelvebf
  \setbox\strutbox\hbox{\vrule height14.6\p@ depth6\p@ width\z@}%
  \setbox\strutbox@\hbox{\lower.86\normallineskiplimit\vbox{%
        \kern-\normallineskiplimit\copy\strutbox}}%
 \setbox\z@\vbox{\hbox{$($}\kern\z@}\bigsize@=2\ht\z@
 \fi
 \normalbaselines\rm\ex@.2326ex\jot5.2\ex@\the\seventeenpoint@}

\catcode`\@=13

\def\({\left(}
\def\){\right)}
\def\[{\left[}
\def\]{\right]}
\def\Gr{\operatorname{Gr}}
\def\sgn{\operatorname{sgn}}
\def\P{\Cal P}

\topmatter 
\title On multiplicities of points on Schubert varieties in
Grassmannians
\endtitle 
\author C.~Krattenthaler\footnote"$^\dagger$"{\hbox{Partially supported by the Austrian
Science Foundation FWF, grant P13190-MAT.}}
\endauthor 
\affil 
Institut f\"ur Mathematik der Universit\"at Wien,\\
Strudlhofgasse 4, A-1090 Wien, Austria.\\
e-mail: KRATT\@Ap.Univie.Ac.At\\
WWW: \tt http://www.mat.univie.ac.at/People/kratt
\endaffil
\address Institut f\"ur Mathematik der Universit\"at Wien,
Strudlhofgasse 4, A-1090 Wien, Austria.
\endaddress
\abstract 
We answer some questions related to multiplicity formulas by Rosenthal
and Zelevinsky and by Lakshmibai and Weyman for points on Schubert
varieties in Grassmannians. In particular, we give combinatorial
interpretations in terms of nonintersecting lattice paths of these
formulas, which makes the equality of the two formulas immediately
obvious. Furthermore we provide an alternative determinantal formula
for these multiplicities, and we show that they count semistandard
tableaux of unusual shapes.
\endabstract
\endtopmatter
\document

\leftheadtext{Christian Krattenthaler}

\subhead 1. Introduction \endsubhead
The {\it multiplicity} of a point on an algebraic variety is an important 
invariant that ``measures" singularity of the point.
Recently, Rosenthal and Zelevinsky \cite{\RoZeAA} gave a
determinantal formula for the multiplicity of a point on
a Schubert variety in a Grassmannian (see Theorem~1). This formula immediately raised
three questions (and, indeed, they are asked in \cite{\RoZeAA,
Remark~5, paragraph after Theorem~1, Remark~7}): 

\roster
\item Is there a direct way to see that the
formula yields positive integers (that is, leaving aside the fact
that the Rosenthal--Zelevinsky
theorem says that it gives multiplicities of singular
points)? 
\item The formula is in form of a binomial determinant. Such
determinants are very common in combinatorics. Is there a
combinatorial interpretation? 
\item Lakshmibai and Weyman \cite{\LaWeAB, Theorem~5.4} 
give a different determinantal formula in a
special case (see Theorem~2). It is not immediately clear why it agrees with the
formula by Rosenthal and Zelevinsky. Is there a straightforward explanation?
\endroster

The purpose of this note is to answer these questions. In reply to
Question~1 we show that, by means of the Lindstr\"om--Gessel--Viennot
theorem \cite{\LindAA, Lemma~1}, \cite{\GeViAB, Theorem~1} (see
Theorem~3),  
the formula of Rosenthal and Zelevinsky
counts certain families of nonintersecting lattice paths. Clearly,
this immediately explains why the formula yields positive integers.
At the same time, this also provides a first answer to Question~2. 
In the special case considered by Lakshmibai and Weyman, we apply an
easy combinatorial transformation to the
families of nonintersecting lattice paths corresponding to the
Rosenthal--Zelevinsky determinant formula and thus convert them into
other families of nonintersecting lattice paths. These latter
families of nonintersecting lattice paths directly yield the
Lakshmibai--Weyman formula, again by means of
the Lindstr\"om--Gessel--Viennot theorem. This answers Question~3.

In addition, we use the ``dual path" idea by Gessel and Viennot
\cite{\GeViAA, Sec.~4} to derive an alternative determinantal formula for the
multiplicities in the general case. As a bonus, this enables us to also
find a combinatorial description of the multiplicities as the numbers of
semistandard tableaux of unusual shapes, thus providing another answer to
Question~2.

In the next section we review the basic definitions and the formulas
by Rosenthal and Zelevinsky and by Lakshmibai and Weyman. Then, in
Section~3, we first recall the Lindstr\"om--Gessel--Viennot theorem
on nonintersecting paths, and then 
explain how to interpret the formulas by 
Rosenthal and Zelevinsky and by Lakshmibai and Weyman in terms of
nonintersecting lattice paths, and why this immediately explains that
they are equivalent in the relevant special case. Finally, in
Section~4, we derive the alternative determinantal formula for the
multiplicities (see Theorem~5) and its interpretation in terms
of semistandard tableaux (see Corollary~6).

\subhead 2. The multiplicity formulas by Rosenthal and Zelevinsky and 
by Lakshmibai and Weyman \endsubhead
Let $d$ and $n$ be positive integers with $0\le d\le n$. The {\it
Grassmannian} $\Gr_d(V)$ is the variety of all $d$-dimensional
subspaces in an $n$-dimensional vector space $V$ (over some
algebraically closed field of arbitrary characteristic). Given an
integer vector $\bold i=(i_1,i_2,\dots,i_d)$, $1\le
i_1<i_2<\dots<i_d\le n$ and a complete flag $\{0\}=V_0\subset
V_1\subset \dots\subset V_n=V$, the {\it Schubert variety} $X_{\bold i}$ is defined by
$$X_{\bold i}=\{W\in \Gr_d(V):\dim(W\cap V_{i_k})\ge k\text { for
}k=1,2,\dots,d\}.$$
The {\it Schubert cell} $X_{\bold i}^\circ$ is an open subset in
$X_{\bold i}$ given by
$$X_{\bold i}^\circ=\{W\in X_{\bold i}:\dim(W\cap V_{i_k-1})= k-1\text { for
}k=1,2,\dots,d\}.$$
It is well-known (see e.g\. 
\cite{\FultAC, Sec.~9.4}) that the Schubert variety $X_{\bold i}$ is
the disjoint union of Schubert cells $X_{\bold j}^\circ$ over all $\bold
j\le \bold i$ (the latter inequality meaning $j_k\le i_k$ for
$k=1,2,\dots,d$). The multiplicity of a point $x$ in
$X_{\bold i}$ is constant on each Schubert cell $X_{\bold
j}^\circ\subset X_{\bold i}$. Following \cite{\RoZeAA} we denote this
multiplicity by $M_{\bold j}(\bold i)$.

The determinantal formula by Rosenthal and Zelevinsky for the
multiplicity $M_{\bold j}(\bold i)$ is the following.
\proclaim{Theorem 1} {\rm(\cite{\RoZeAA, Theorem~1})} The multiplicity 
$M_{\bold j}(\bold i)$ of a point $x\in X_{\bold j}^\circ\subset 
X_{\bold i}$ is given by
$$M_{\bold j}(\bold i)=(-1)^{s_1+\dots+s_d}\det_{1\le p,q\le d}
\left(\binom {i_q}{p-1-s_q}\right),\tag1$$
where $s_q=\vert \{\ell:i_q<j_\ell\}\vert$. 
\endproclaim

In the special case that $\bold j=(1,2,\dots,d)$, 
Lakshmibai and Weyman have given a different determinant formula.
In the statement of their theorem we use standard partition
terminology (see e.g\. \cite{\MacdAC, Ch.~I, Sec.~1}).
\proclaim{Theorem 2} {\rm(\cite{\LaWeAB, Theorem~5.4})}
The multiplicity 
$M_{(1,2,\dots,d)}(\bold i)$ of a point $x\in X_{(1,2,\dots,d)}^\circ\subset 
X_{\bold i}$ is given by
$$M_{(1,2,\dots,d)}(\bold i)=\det_{1\le p,q\le r}\left(\binom
{\alpha_p+\beta_q}{\alpha_p}\right),\tag2$$
where $(\alpha_1,\dots,\alpha_r\mid \beta_1,\dots,\beta_r)$ is the Frobenius
notation of the partition $\lambda=(i_d-d,\dots,i_2-2,i_1-1)$. 
\endproclaim

\subhead 3. Multiplicities count nonintersecting lattice paths \endsubhead
We start by recalling the main theorem on nonintersecting lattice
paths, due to Lindstr\"om, and Gessel and
Viennot. 
\proclaim{Theorem 3} {\rm (\cite{\LindAA, Lemma~1}, \cite{\GeViAB,
Theorem~1})}
Let $G$ be any acyclic directed graph. 
Let $A_1,A_2,\dots, A_d, E_1, E_2,\dots , E_d$ be vertices of $G$.
Then, with $S_d$ denoting the group of permutations of $\{1,2,\dots,d\}$, 
the following identity holds:
$$
\det_{1\le p,q \le d}{\(\P (A_p \to E_q)\)}=\sum _{\sigma\in S_d} ^{}{
(\sgn \sigma )\cdot \P ({\bold A} \to {\bold E_\sigma}, \text {\rm
nonint.})},
\tag3$$
where $\P(A\to E)$ denotes the number of paths from $A$ to $E$ in
$G$, and where $\P ({\bold A} \to {\bold E_\sigma}, \text {\rm
nonint.})$ denotes the number of all families $(P_1,P_2,\dots,P_d)$
of paths in $G$, $P_\ell$ running from $A_\ell$ to $E_{\sigma(\ell)}$,
$\ell=1,2,\dots,d$, which are
{\rm nonintersecting}. A family of paths is called 
nonintersecting if no two paths
in the family have a point in common.
\endproclaim

The most commonly used instance of this rather general theorem arises when
the starting and end points are in a position such that the only
nonvanishing term on the right-hand side of (3) is the one for
$\sigma$ equal to the identity permutation.

\proclaim{Corollary 4} {\rm (\cite{\GeViAB,
Cor.~2})} In addition to the assumptions in Theorem~3,
assume that for any $i<j$ and $k<l$ any path from $A_i$ to $E_l$
intersects any path from $A_j$ to $A_k$. Then the number of 
all families $(P_1,P_2,\dots,P_d)$
of nonintersecting paths in $G$, $P_\ell$ running from $A_\ell$ to
$E_{\sigma(\ell)}$,
$\ell=1,2,\dots,d$, is equal to
$$\det_{1\le p,q \le d}{\(\P (A_p \to E_q)\)}.$$
\endproclaim

In view of the above corollary, the determinant (2) by Lakshmibai and Weyman
has an obvious interpretation in terms of
nonintersecting lattice paths consisting of horizontal and vertical
steps in the positive direction: It counts the number of all families
$(P_1,P_2,\dots,P_r)$ of nonintersecting lattice paths, the path $P_\ell$
running from $(-\beta_\ell,0)$ to $(0,\alpha_\ell)$, $\ell=1,2,\dots,r$.
See Figure~1 for an example, where $d=7$,
$(i_1,i_2,\dots,i_7)=(3,5,9,10,14,15,17)$, and, hence, 
$\lambda=(10,9,9,6,6,3,2)=(9,7,6,2,1\mid 6,5,3,1,0)$.

\midinsert
\vskip8pt
\vbox{
$$
\Gitter(8,10)(0,0)
\PfadDicke{.5pt}
\Pfad(7,-1),22222222222\endPfad
\Pfad(-1,0),1111111111\endPfad
\PfadDicke{1.2pt}
\Pfad(1,0),222222112211211\endPfad
\Pfad(2,0),222121222111\endPfad
\Pfad(4,0),212221212\endPfad
\Pfad(6,0),221\endPfad
\Pfad(7,0),2\endPfad
\DickPunkt(1,0)
\DickPunkt(2,0)
\DickPunkt(4,0)
\DickPunkt(6,0)
\DickPunkt(7,0)
\DickPunkt(7,1)
\DickPunkt(7,2)
\DickPunkt(7,6)
\DickPunkt(7,7)
\DickPunkt(7,9)
\hskip3.5cm
$$
\centerline{\eightpoint Figure 1}
}
\vskip8pt
\endinsert

The formula by Rosenthal and Zelevinsky can also be interpreted in
terms of nonintersecting lattice paths. 
By Theorem~3, the determinant (1) counts the weighted sum of all families
$(Q_1,Q_2,\dots,Q_d)$ of nonintersecting lattice paths, where the
path $Q_\ell$ runs from $(-\ell+1,\ell-1)$,
to $(-s_{\sigma(\ell)},s_{\sigma(\ell)}+i_{\sigma(\ell)})$, 
$\ell=1,2,\dots,d$, for some permutation $\sigma\in S_d$, 
and where the weight of a path family is defined
as the sign of $\sigma$.
See Figure~2 for an example with the same parameters, 
in which $s_1=4$, $s_2=2$, and
$s_3=\dots=s_7=0$. 
However, it is not difficult to see that the permutation $\sigma$ is
in fact always the same, and that the sign of $\sigma$ is
$(-1)^{s_1+\dots+s_d}$.

\midinsert
\vskip8pt
\vbox{
$$
\Gitter(8,18)(0,0)
\PfadDicke{.5pt}
\Pfad(7,-1),222222222222222222\endPfad
\Pfad(-1,0),1111111111\endPfad
\PfadDicke{1.2pt}
\Pfad(1,6),22222222112211211\endPfad
\Pfad(2,5),222222121222111\endPfad
\Pfad(3,4),222\endPfad
\Pfad(4,3),22222212221212\endPfad
\Pfad(5,2),22222\endPfad
\Pfad(6,1),2222222221\endPfad
\Pfad(7,0),222222222\endPfad
\SPfad(-1,8),111111111\endSPfad
\DickPunkt(1,6)
\DickPunkt(2,5)
\DickPunkt(3,4)
\DickPunkt(4,3)
\DickPunkt(5,2)
\DickPunkt(6,1)
\DickPunkt(3,7)
\DickPunkt(5,7)
\DickPunkt(7,0)
\DickPunkt(7,9)
\DickPunkt(7,10)
\DickPunkt(7,14)
\DickPunkt(7,15)
\DickPunkt(7,17)
\hskip3.5cm
$$
\centerline{\eightpoint Figure 2}
}
\vskip8pt
\endinsert

This gives a combinatorial interpretation of the
Rosenthal--Zelevinsky formula for
{\it any} $\bold i$ and $\bold j$, namely as the number
of all families
$(Q_1,Q_2,\dots,Q_d)$ of nonintersecting lattice paths, where the
path $Q_\ell$ runs from $(-\ell+1,\ell-1)$,
to $(-s_{\sigma(\ell)},s_{\sigma(\ell)}+i_{\sigma(\ell)})$, 
$\ell=1,2,\dots,d$, for some permutation $\sigma\in S_d$ (which is
uniquely determined).

In the special case, however, that
$\bold j=(1,2,\dots,d)$, we have $s_\ell=d-i_\ell$ as long as
$i_\ell\le d$, so that the end points of the paths are either
$(i_\ell-d,d)$ or $(0,i_\ell)$. 
It is now easy to bijectively map these path families to the former.
The path $Q_1$ starts at $(0,0)$. If the path family should be
nonintersecting, then the
only possibility for $Q_1$ is to run from $(0,0)$ to $(0,i_m)$, where
$m$ is minimal such that $s_m=0$. (This minimum exists since we must
have $s_d=0$.) In particular, since we must have $i_m>d$, the path
$Q_1$ starts with $d+1$ vertical steps. This forces all the other
paths (if they want to be nonintersecting) to also be vertical until
they reach height $d+1$. If they do not reach height $d+1$, then
they terminate in
a point $(i_\ell-d,d)$. All this is clearly visible in Figure~2. 

Hence, we may without loss of any information cut off the path
portions until height $d+1$, respectively remove the paths that even
do not reach that height. (The dotted line in Figure~2 indicates the
line of height $d+1=8$ along which the cut is performed.) 
What remains is a family of paths that is
familiar from the formula of Lakshmibai and Weyman, see Figure~1.

This makes it obvious why formula (1)
with $\bold j=(1,2,\dots,d)$ and formula (2) agree.

\subhead 4. Dual paths and multiplicities as numbers of semistandard
tableaux \endsubhead
In this section we derive an alternative determinantal formula for
the multiplicities, and we show that they count certain semistandard
tableaux.

If we would try to use the ideas of the previous section in the
general case (i.e., for arbitrary $\bold j\le \bold i$), then we
first discover that, in general, it is not true that $s_\ell+i_\ell$ is
equal to $d$ or $i_\ell$. We can just say that
the numbers $s_\ell+i_\ell$, $\ell=1,2,\dots,d$, are
weakly increasing. Yet, we try the same construction.

Recall from the previous section that
the multiplicity $M_{\bold j}(\bold i)$ is equal to
the number of all families
$(Q_1,Q_2,\dots,Q_d)$ of nonintersecting lattice paths, where the
path $Q_\ell$ runs from $(-\ell+1,\ell-1)$,
to $(-s_{\sigma(\ell)},s_{\sigma(\ell)}+i_{\sigma(\ell)})$, 
$\ell=1,2,\dots,d$, for some permutation $\sigma\in S_d$. (As we
remarked, there is a unique permutation $\sigma$ for which such families
of nonintersecting lattice paths exist).

Let us consider an example: $d=9$, $\bold i=(4,6,7,13,14,17,19,20,21)$
and $\bold j=(1,2,4,7,10,12,13,15,16)$. Then $\bold
s=(6,6,5,2,2,0,0,0,0)$. 
Figure~3 shows a typical family for this choice of $\bold i$ and
$\bold j$.

\midinsert
\vskip8pt
\vbox{
$$
\Gitter(10,22)(0,0)
\PfadDicke{.5pt}
\Pfad(9,-1),2222222222222222222222\endPfad
\Pfad(-1,0),11111111111\endPfad
\PfadDicke{1.2pt}
\Pfad(1,8),222221222112212221111\endPfad
\Pfad(2,7),222212\endPfad
\Pfad(3,6),2222\endPfad
\Pfad(4,5),2222222\endPfad
\Pfad(5,4),22222222222221221211\endPfad
\Pfad(6,3),22222222222221\endPfad
\Pfad(7,2),2222222222222\endPfad
\Pfad(8,1),2222222222222222212\endPfad
\Pfad(9,0),22222222222222222\endPfad
\SPfad(1,10),112212221112211\endSPfad
\DickPunkt(1,8)
\DickPunkt(2,7)
\DickPunkt(3,6)
\DickPunkt(4,5)
\DickPunkt(5,4)
\DickPunkt(6,3)
\DickPunkt(7,2)
\DickPunkt(8,1)
\DickPunkt(9,0)
\DickPunkt(3,10)
\DickPunkt(3,12)
\DickPunkt(4,12)
\DickPunkt(7,15)
\DickPunkt(7,16)
\DickPunkt(9,17)
\DickPunkt(9,19)
\DickPunkt(9,20)
\DickPunkt(9,21)
\hskip3.5cm
$$
\centerline{\eightpoint Figure 3}
}
\vskip8pt
\endinsert

Again, large parts of the initial vertical portions of the paths are
forced. To be precise, let $\{s_1,s_2,\dots,s_d\}=\{v_1,v_2,\dots,v_r\}$, 
$v_1<v_2<\dots<v_r$. (I.e., the $v_\ell$'s are the distinct values
that are attained by the $s_\ell$'s.) For each $\ell$,
consider the
bottom-most end point $(-s_q,s_q+i_q)$ with $s_q=v_\ell$. Then we may
cut off the vertical portions below height $s_q+i_q$ of the paths
with starting points $(-v_{\ell+1}+1,v_{\ell+1}-1)$, \dots, 
$(-v_{\ell}-1,v_{\ell}+1)$, 
$(-v_{\ell},v_{\ell})$. These cuts are indicated by the dotted line
in Figure~3. (I.e., portions below dotted lines can be
omitted.) The result, after the cuts, is shown in Figure~4.

\def\GPunkt(#1,#2){\unskip
  \raise#2 \Einheit\hbox to0pt{\hskip#1 \Einheit
          \raise-2.5pt\hbox to0pt{\hss$\ssize\bullet$\hss}\hss}}

\midinsert
\vskip8pt
\vbox{
$$
\PfadDicke{.5pt}
\Pfad(9,-1),2222222222222\endPfad
\Pfad(-1,0),11111111111\endPfad
\PfadDicke{1.2pt}
\Pfad(1,0),2221222112212221111\endPfad
\Pfad(2,0),212\endPfad
\Pfad(5,5),221221211\endPfad
\Pfad(6,5),21\endPfad
\Pfad(8,7),212\endPfad
\SPfad(1,0),112212221112211\endSPfad
\DickPunkt(1,0)
\DickPunkt(2,0)
\DickPunkt(5,5)
\DickPunkt(6,5)
\DickPunkt(8,7)
\DickPunkt(3,0)
\DickPunkt(3,2)
\DickPunkt(4,2)
\DickPunkt(7,5)
\DickPunkt(7,6)
\DickPunkt(9,7)
\DickPunkt(9,9)
\DickPunkt(9,10)
\DickPunkt(9,11)
\GPunkt(1,1)
\GPunkt(2,1)
\GPunkt(3,1)
\GPunkt(1,2)
\GPunkt(2,2)
\GPunkt(1,3)
\GPunkt(2,3)
\GPunkt(3,3)
\GPunkt(4,3)
\GPunkt(1,4)
\GPunkt(2,4)
\GPunkt(3,4)
\GPunkt(4,4)
\GPunkt(1,5)
\GPunkt(2,5)
\GPunkt(3,5)
\GPunkt(4,5)
\GPunkt(1,6)
\GPunkt(2,6)
\GPunkt(3,6)
\GPunkt(4,6)
\GPunkt(5,6)
\GPunkt(6,6)
\GPunkt(1,7)
\GPunkt(2,7)
\GPunkt(3,7)
\GPunkt(4,7)
\GPunkt(5,7)
\GPunkt(6,7)
\GPunkt(7,7)
\GPunkt(8,7)
\GPunkt(1,8)
\GPunkt(2,8)
\GPunkt(3,8)
\GPunkt(4,8)
\GPunkt(5,8)
\GPunkt(6,8)
\GPunkt(7,8)
\GPunkt(8,8)
\GPunkt(9,8)
\GPunkt(1,9)
\GPunkt(2,9)
\GPunkt(3,9)
\GPunkt(4,9)
\GPunkt(5,9)
\GPunkt(6,9)
\GPunkt(7,9)
\GPunkt(8,9)
\GPunkt(9,9)
\GPunkt(1,10)
\GPunkt(2,10)
\GPunkt(3,10)
\GPunkt(4,10)
\GPunkt(5,10)
\GPunkt(6,10)
\GPunkt(7,10)
\GPunkt(8,10)
\GPunkt(9,10)
\GPunkt(1,11)
\GPunkt(2,11)
\GPunkt(3,11)
\GPunkt(4,11)
\GPunkt(5,11)
\GPunkt(6,11)
\GPunkt(7,11)
\GPunkt(8,11)
\GPunkt(9,11)
\hskip3.5cm
$$
\centerline{\eightpoint Figure 4}
}
\vskip8pt
\endinsert

Now we may write down a Lindstr\"om--Gessel--Viennot determinant for
these (new) starting points and (old) end points. The result would
again be a determinant of binomials. However, it seems
that it would require a considerable 
amount of notation to explicitely express
what the new starting points are. 

Moreover, this would not fully
correspond to the Lakshmibai-Weyman formula because, 
to obtain the Lakshmibai-Weyman formula, one had to also drop the
paths of length zero. If we would do that in the above picture, then
we would have to restrict the paths explicitly to the indicated
ladder-shaped region. If we would not do that then, after removal of the
starting and end points corresponding to the zero length paths (which
opens ``holes"), there are now more possibilities to connect the
(reduced set of) starting points to the (reduced set of) end points
by nonintersecting lattice paths, and the corresponding permutation
$\sigma$ would not be unique anymore. So, we could still write down a 
Lindstr\"om--Gessel--Viennot determinant (that would now correspond to
the Lakshibai-Weyman formula), however the entries would not be
binomials anymore, they would count paths with given starting and end
points {\it that stay in this ladder-shaped region}, for which no closed
formula is available.

\midinsert
\vskip8pt
\vbox{
$$
\PfadDicke{.5pt}
\Pfad(9,-1),2222222222222\endPfad
\PfadDicke{1.2pt}
\Pfad(1,0),2221222112212221111\endPfad
\Pfad(2,0),212\endPfad
\Pfad(5,5),221221211\endPfad
\Pfad(6,5),21\endPfad
\Pfad(8,7),212\endPfad
\SPfad(1,11),66666666464\endSPfad
\SPfad(2,11),666664666\endSPfad
\SPfad(3,11),666664666\endSPfad
\SPfad(4,11),666444\endSPfad
\SPfad(5,11),46466\endSPfad
\SPfad(6,11),4464\endSPfad
\SPfad(7,11),44\endSPfad
\SPfad(8,11),4\endSPfad
\Kreis(1,11)
\Kreis(2,11)
\Kreis(3,11)
\Kreis(4,11)
\Kreis(5,11)
\Kreis(6,11)
\Kreis(7,11)
\Kreis(8,11)
\Kreis(9,11)
\DickPunkt(1,0)
\DickPunkt(2,0)
\DickPunkt(5,5)
\DickPunkt(6,5)
\DickPunkt(8,7)
\DickPunkt(3,0)
\DickPunkt(3,2)
\DickPunkt(4,2)
\DickPunkt(7,5)
\DickPunkt(7,6)
\DickPunkt(9,7)
\DickPunkt(9,9)
\DickPunkt(9,10)
\GPunkt(1,1)
\GPunkt(2,1)
\GPunkt(3,1)
\GPunkt(1,2)
\GPunkt(2,2)
\GPunkt(1,3)
\GPunkt(2,3)
\GPunkt(3,3)
\GPunkt(4,3)
\GPunkt(1,4)
\GPunkt(2,4)
\GPunkt(3,4)
\GPunkt(4,4)
\GPunkt(1,5)
\GPunkt(2,5)
\GPunkt(3,5)
\GPunkt(4,5)
\GPunkt(1,6)
\GPunkt(2,6)
\GPunkt(3,6)
\GPunkt(4,6)
\GPunkt(5,6)
\GPunkt(6,6)
\GPunkt(1,7)
\GPunkt(2,7)
\GPunkt(3,7)
\GPunkt(4,7)
\GPunkt(5,7)
\GPunkt(6,7)
\GPunkt(7,7)
\GPunkt(8,7)
\GPunkt(1,8)
\GPunkt(2,8)
\GPunkt(3,8)
\GPunkt(4,8)
\GPunkt(5,8)
\GPunkt(6,8)
\GPunkt(7,8)
\GPunkt(8,8)
\GPunkt(9,8)
\GPunkt(1,9)
\GPunkt(2,9)
\GPunkt(3,9)
\GPunkt(4,9)
\GPunkt(5,9)
\GPunkt(6,9)
\GPunkt(7,9)
\GPunkt(8,9)
\GPunkt(9,9)
\GPunkt(1,10)
\GPunkt(2,10)
\GPunkt(3,10)
\GPunkt(4,10)
\GPunkt(5,10)
\GPunkt(6,10)
\GPunkt(7,10)
\GPunkt(8,10)
\GPunkt(9,10)
\GPunkt(1,11)
\GPunkt(2,11)
\GPunkt(3,11)
\GPunkt(4,11)
\GPunkt(5,11)
\GPunkt(6,11)
\GPunkt(7,11)
\GPunkt(8,11)
\GPunkt(9,11)
\hbox{\hskip7cm}
\PfadDicke{.5pt}
\Pfad(9,-1),2222222222222\endPfad
\PfadDicke{1.2pt}
\SPfad(1,0),2221222112212221111\endSPfad
\SPfad(2,0),212\endSPfad
\SPfad(5,5),221221211\endSPfad
\SPfad(6,5),21\endSPfad
\SPfad(8,7),212\endSPfad
\Pfad(1,11),66666666464\endPfad
\Pfad(2,11),666664666\endPfad
\Pfad(3,11),666664666\endPfad
\Pfad(4,11),666444\endPfad
\Pfad(5,11),46466\endPfad
\Pfad(6,11),4464\endPfad
\Pfad(7,11),44\endPfad
\Pfad(8,11),4\endPfad
\Kreis(1,11)
\Kreis(2,11)
\Kreis(3,11)
\Kreis(4,11)
\Kreis(5,11)
\Kreis(6,11)
\Kreis(7,11)
\Kreis(8,11)
\Kreis(9,11)
\DickPunkt(1,0)
\DickPunkt(2,0)
\DickPunkt(5,5)
\DickPunkt(6,5)
\DickPunkt(8,7)
\DickPunkt(3,0)
\DickPunkt(3,2)
\DickPunkt(4,2)
\DickPunkt(7,5)
\DickPunkt(7,6)
\DickPunkt(9,7)
\DickPunkt(9,9)
\DickPunkt(9,10)
\GPunkt(1,1)
\GPunkt(2,1)
\GPunkt(3,1)
\GPunkt(1,2)
\GPunkt(2,2)
\GPunkt(1,3)
\GPunkt(2,3)
\GPunkt(3,3)
\GPunkt(4,3)
\GPunkt(1,4)
\GPunkt(2,4)
\GPunkt(3,4)
\GPunkt(4,4)
\GPunkt(1,5)
\GPunkt(2,5)
\GPunkt(3,5)
\GPunkt(4,5)
\GPunkt(1,6)
\GPunkt(2,6)
\GPunkt(3,6)
\GPunkt(4,6)
\GPunkt(5,6)
\GPunkt(6,6)
\GPunkt(1,7)
\GPunkt(2,7)
\GPunkt(3,7)
\GPunkt(4,7)
\GPunkt(5,7)
\GPunkt(6,7)
\GPunkt(7,7)
\GPunkt(8,7)
\GPunkt(1,8)
\GPunkt(2,8)
\GPunkt(3,8)
\GPunkt(4,8)
\GPunkt(5,8)
\GPunkt(6,8)
\GPunkt(7,8)
\GPunkt(8,8)
\GPunkt(9,8)
\GPunkt(1,9)
\GPunkt(2,9)
\GPunkt(3,9)
\GPunkt(4,9)
\GPunkt(5,9)
\GPunkt(6,9)
\GPunkt(7,9)
\GPunkt(8,9)
\GPunkt(9,9)
\GPunkt(1,10)
\GPunkt(2,10)
\GPunkt(3,10)
\GPunkt(4,10)
\GPunkt(5,10)
\GPunkt(6,10)
\GPunkt(7,10)
\GPunkt(8,10)
\GPunkt(9,10)
\GPunkt(1,11)
\GPunkt(2,11)
\GPunkt(3,11)
\GPunkt(4,11)
\GPunkt(5,11)
\GPunkt(6,11)
\GPunkt(7,11)
\GPunkt(8,11)
\GPunkt(9,11)
\hskip3.5cm
$$
\centerline{\eightpoint Figure 5}
}
\vskip8pt
\endinsert

On the other hand, as announced, we may now
introduce what Gessel and Viennot \cite{\GeViAA, Sec.~4} call ``dual
paths," see Figure~5. We mark starting points $(d-\ell,i_d)$,
$\ell=1,2,\dots,d$. (They are indicated by circles in Figure~5.) 
Then, for each of these points, we move
vertically downwards, unless we hit one of the existing paths. If the
latter happens, then we continue by a diagonal step $(1,-1)$, etc. It is
not difficult to see that in that manner we connect $(d-\ell,i_\ell)$
with $(-s_\ell,s_\ell+i_\ell)$, $\ell=1,2,\dots,d$. The
resulting paths in our running example are indicated by dotted lines 
in the left picture in Figure~5. In the right picture I have just
interchanged the roles of the two families of paths.

Now we deform the lattice slightly, so that the newly introduced paths
become orthogonal paths, see Figure~6.

\midinsert
\vskip8pt
\vbox{
$$
\PfadDicke{.5pt}
\PfadDicke{1.2pt}
\SPfad(1,0),2223222332232223333\endSPfad
\SPfad(2,1),232\endSPfad
\SPfad(5,9),223223233\endSPfad
\SPfad(6,10),23\endSPfad
\SPfad(8,14),232\endSPfad
\Pfad(1,11),66666666161\endPfad
\Pfad(2,12),666661666\endPfad
\Pfad(3,13),666661666\endPfad
\Pfad(4,14),666111\endPfad
\Pfad(5,15),16166\endPfad
\Pfad(6,16),1161\endPfad
\Pfad(7,17),11\endPfad
\Pfad(8,18),1\endPfad
\Kreis(1,11)
\Kreis(2,12)
\Kreis(3,13)
\Kreis(4,14)
\Kreis(5,15)
\Kreis(6,16)
\Kreis(7,17)
\Kreis(8,18)
\Kreis(9,19)
\DickPunkt(1,0)
\DickPunkt(2,1)
\DickPunkt(5,9)
\DickPunkt(6,10)
\DickPunkt(8,14)
\DickPunkt(3,2)
\DickPunkt(3,4)
\DickPunkt(4,5)
\DickPunkt(7,11)
\DickPunkt(7,12)
\DickPunkt(9,15)
\DickPunkt(9,17)
\DickPunkt(9,18)
\GPunkt(1,0)
\GPunkt(1,1)
\GPunkt(1,2)
\GPunkt(1,3)
\GPunkt(1,4)
\GPunkt(1,5)
\GPunkt(1,6)
\GPunkt(1,7)
\GPunkt(1,8)
\GPunkt(1,9)
\GPunkt(1,10)
\GPunkt(1,11)
\GPunkt(2,1)
\GPunkt(2,2)
\GPunkt(2,3)
\GPunkt(2,4)
\GPunkt(2,5)
\GPunkt(2,6)
\GPunkt(2,7)
\GPunkt(2,8)
\GPunkt(2,9)
\GPunkt(2,10)
\GPunkt(2,11)
\GPunkt(2,12)
\GPunkt(3,2)
\GPunkt(3,3)
\GPunkt(3,4)
\GPunkt(3,5)
\GPunkt(3,6)
\GPunkt(3,7)
\GPunkt(3,8)
\GPunkt(3,9)
\GPunkt(3,10)
\GPunkt(3,11)
\GPunkt(3,12)
\GPunkt(3,13)
\GPunkt(4,5)
\GPunkt(4,6)
\GPunkt(4,7)
\GPunkt(4,8)
\GPunkt(4,9)
\GPunkt(4,10)
\GPunkt(4,11)
\GPunkt(4,12)
\GPunkt(4,13)
\GPunkt(4,14)
\GPunkt(5,9)
\GPunkt(5,10)
\GPunkt(5,11)
\GPunkt(5,12)
\GPunkt(5,13)
\GPunkt(5,14)
\GPunkt(5,15)
\GPunkt(6,10)
\GPunkt(6,11)
\GPunkt(6,12)
\GPunkt(6,13)
\GPunkt(6,14)
\GPunkt(6,15)
\GPunkt(6,16)
\GPunkt(7,11)
\GPunkt(7,12)
\GPunkt(7,13)
\GPunkt(7,14)
\GPunkt(7,15)
\GPunkt(7,16)
\GPunkt(7,17)
\GPunkt(8,14)
\GPunkt(8,15)
\GPunkt(8,16)
\GPunkt(8,17)
\GPunkt(8,18)
\GPunkt(9,15)
\GPunkt(9,16)
\GPunkt(9,17)
\GPunkt(9,18)
\GPunkt(9,19)
\hbox{\hskip7cm}
\PfadDicke{.5pt}
\PfadDicke{1.2pt}
\Pfad(1,11),66666666161\endPfad
\Pfad(2,12),666661666\endPfad
\Pfad(3,13),666661666\endPfad
\Pfad(4,14),666111\endPfad
\Pfad(5,15),16166\endPfad
\Pfad(6,16),1161\endPfad
\Pfad(7,17),11\endPfad
\Pfad(8,18),1\endPfad
\Kreis(1,11)
\Kreis(2,12)
\Kreis(3,13)
\Kreis(4,14)
\Kreis(5,15)
\Kreis(6,16)
\Kreis(7,17)
\Kreis(8,18)
\Kreis(9,19)
\DickPunkt(3,2)
\DickPunkt(3,4)
\DickPunkt(4,5)
\DickPunkt(7,11)
\DickPunkt(7,12)
\DickPunkt(9,15)
\DickPunkt(9,17)
\DickPunkt(9,18)
\GPunkt(1,0)
\GPunkt(1,1)
\GPunkt(1,2)
\GPunkt(1,3)
\GPunkt(1,4)
\GPunkt(1,5)
\GPunkt(1,6)
\GPunkt(1,7)
\GPunkt(1,8)
\GPunkt(1,9)
\GPunkt(1,10)
\GPunkt(1,11)
\GPunkt(2,1)
\GPunkt(2,2)
\GPunkt(2,3)
\GPunkt(2,4)
\GPunkt(2,5)
\GPunkt(2,6)
\GPunkt(2,7)
\GPunkt(2,8)
\GPunkt(2,9)
\GPunkt(2,10)
\GPunkt(2,11)
\GPunkt(2,12)
\GPunkt(3,2)
\GPunkt(3,3)
\GPunkt(3,4)
\GPunkt(3,5)
\GPunkt(3,6)
\GPunkt(3,7)
\GPunkt(3,8)
\GPunkt(3,9)
\GPunkt(3,10)
\GPunkt(3,11)
\GPunkt(3,12)
\GPunkt(3,13)
\GPunkt(4,5)
\GPunkt(4,6)
\GPunkt(4,7)
\GPunkt(4,8)
\GPunkt(4,9)
\GPunkt(4,10)
\GPunkt(4,11)
\GPunkt(4,12)
\GPunkt(4,13)
\GPunkt(4,14)
\GPunkt(5,9)
\GPunkt(5,10)
\GPunkt(5,11)
\GPunkt(5,12)
\GPunkt(5,13)
\GPunkt(5,14)
\GPunkt(5,15)
\GPunkt(6,10)
\GPunkt(6,11)
\GPunkt(6,12)
\GPunkt(6,13)
\GPunkt(6,14)
\GPunkt(6,15)
\GPunkt(6,16)
\GPunkt(7,11)
\GPunkt(7,12)
\GPunkt(7,13)
\GPunkt(7,14)
\GPunkt(7,15)
\GPunkt(7,16)
\GPunkt(7,17)
\GPunkt(8,14)
\GPunkt(8,15)
\GPunkt(8,16)
\GPunkt(8,17)
\GPunkt(8,18)
\GPunkt(9,15)
\GPunkt(9,16)
\GPunkt(9,17)
\GPunkt(9,18)
\GPunkt(9,19)
\hskip3.5cm
$$
\centerline{\eightpoint Figure 6}
}
\vskip8pt
\endinsert

So, what we obtain finally 
is a family $(R_1,R_2,\dots,R_d)$ of nonintersecting lattice paths, where the
path $R_\ell$ is a path consisting of horizontal unit steps in the
positive direction
and vertical unit steps in the negative direction 
and runs from $(-d+\ell,i_d+\ell)$,
to $(-s_\ell,i_\ell+d)$, $\ell=1,2,\dots,d$. Hence, again by 
Corollary~4, we obtain that the number of
these families of nonintersecting lattice paths is equal to the determinant
$$\det_{1\le p,q\le d}\left(\binom {i_d-i_p-s_q}{d-p-s_q}\right).$$
This proves the following alternative to the result by Rosenthal and
Zelevinsky.
\proclaim{Theorem 5} The multiplicity 
$M_{\bold j}(\bold i)$ of a point $x\in X_{\bold j}^\circ\subset 
X_{\bold i}$ is given by
$$\det_{1\le p,q\le d}\left(\binom
{i_d-i_p-s_q}{d-p-s_q}\right).\tag4$$
where, again, $s_q=\vert \{\ell:i_q<j_\ell\}\vert$. 
\endproclaim
\remark{Remark} 
In fact, this determinant could be restricted to $1\le p,q\le d-1$,
because $\binom {i_d-i_d-s_q}{d-d-s_q}=\binom
{-s_q}{-s_q}=\delta_{q,d}$ (with $\delta_{q,d}$ the Kronecker delta). 
This corresponds to the combinatorial
fact that path $P_d$ is a zero length path. See the top-right of
Figure~6.
\endremark

\midinsert
\vskip8pt
\vbox{
$$
\PfadDicke{.5pt}
\PfadDicke{1.2pt}
\Pfad(1,11),66666666161\endPfad
\Pfad(2,12),666661666\endPfad
\Pfad(3,13),666661666\endPfad
\Pfad(4,14),666111\endPfad
\Pfad(5,15),16166\endPfad
\Pfad(6,16),1161\endPfad
\Pfad(7,17),11\endPfad
\Pfad(8,18),1\endPfad
\Kreis(1,11)
\Kreis(2,12)
\Kreis(3,13)
\Kreis(4,14)
\Kreis(5,15)
\Kreis(6,16)
\Kreis(7,17)
\Kreis(8,18)
\Kreis(9,19)
\DickPunkt(3,2)
\DickPunkt(3,4)
\DickPunkt(4,5)
\DickPunkt(7,11)
\DickPunkt(7,12)
\DickPunkt(9,15)
\DickPunkt(9,17)
\DickPunkt(9,18)
\GPunkt(1,0)
\GPunkt(1,1)
\GPunkt(1,2)
\GPunkt(1,3)
\GPunkt(1,4)
\GPunkt(1,5)
\GPunkt(1,6)
\GPunkt(1,7)
\GPunkt(1,8)
\GPunkt(1,9)
\GPunkt(1,10)
\GPunkt(1,11)
\GPunkt(2,1)
\GPunkt(2,2)
\GPunkt(2,3)
\GPunkt(2,4)
\GPunkt(2,5)
\GPunkt(2,6)
\GPunkt(2,7)
\GPunkt(2,8)
\GPunkt(2,9)
\GPunkt(2,10)
\GPunkt(2,11)
\GPunkt(2,12)
\GPunkt(3,2)
\GPunkt(3,3)
\GPunkt(3,4)
\GPunkt(3,5)
\GPunkt(3,6)
\GPunkt(3,7)
\GPunkt(3,8)
\GPunkt(3,9)
\GPunkt(3,10)
\GPunkt(3,11)
\GPunkt(3,12)
\GPunkt(3,13)
\GPunkt(4,5)
\GPunkt(4,6)
\GPunkt(4,7)
\GPunkt(4,8)
\GPunkt(4,9)
\GPunkt(4,10)
\GPunkt(4,11)
\GPunkt(4,12)
\GPunkt(4,13)
\GPunkt(4,14)
\GPunkt(5,9)
\GPunkt(5,10)
\GPunkt(5,11)
\GPunkt(5,12)
\GPunkt(5,13)
\GPunkt(5,14)
\GPunkt(5,15)
\GPunkt(6,10)
\GPunkt(6,11)
\GPunkt(6,12)
\GPunkt(6,13)
\GPunkt(6,14)
\GPunkt(6,15)
\GPunkt(6,16)
\GPunkt(7,11)
\GPunkt(7,12)
\GPunkt(7,13)
\GPunkt(7,14)
\GPunkt(7,15)
\GPunkt(7,16)
\GPunkt(7,17)
\GPunkt(8,14)
\GPunkt(8,15)
\GPunkt(8,16)
\GPunkt(8,17)
\GPunkt(8,18)
\GPunkt(9,15)
\GPunkt(9,16)
\GPunkt(9,17)
\GPunkt(9,18)
\GPunkt(9,19)
\Label\ro{1}(5,15)
\Label\ro{1}(6,16)
\Label\ro{1}(7,17)
\Label\ro{1}(8,18)
\Label\ro{2}(8,17)
\Label\ro{2}(7,16)
\Label\ro{3}(6,14)
\Label\ro{4}(8,15)
\Label\ro{4}(4,11)
\Label\ro{5}(5,11)
\Label\ro{6}(6,11)
\Label\ro{6}(2,7)
\Label\ro{6}(3,8)
\Label\ro{9}(1,3)
\Label\ro{11}(2,2)
\hskip7cm
\raise5cm\hbox{$
\smatrix \format\sa\c\s\c\s\c\s\c\s\c\s\c\s\c\s\c\se\\
\hlinefor{17}\\
&1&&1&&1&&1&&4&&6&&6&&9&\\
\hlinefor{17}\\
\omit&&&2&&2&&3&&5&& &\omit& &&11&\\
\omit& \omit&\hlinefor9& \omit&\omit& \omit&\hlinefor3\\
\omit& &\omit& &&4&& &&6&& \\
\omit& \omit&\omit& \omit&\hlinefor3& \omit &\hlinefor3
\endsmatrix$}
$$
\centerline{\eightpoint Figure 7}
}
\vskip8pt
\endinsert

As a bonus, we are now able to derive a combinatorial interpretation
of the multiplicities $M_{\bold j}(\bold i)$ in terms of
certain semistandard tableaux. If we label horizontal steps along diagonals
by 1, 2, \dots, respectively, as indicated in Figure~7 (i.e., along the
first diagonal of horizontal steps immediately to the right of the
starting points, the latter being indicated by circles, horizontal steps
are labelled by 1, along the next
diagonal horizontal steps are labelled by 2, etc.), and then read the labels along each
path and form columns out of it, then we obtain an array, as shown in
Figure~7 in our example, which has the following properties:
\roster
\item The length of column $\ell$ is $\ell-s_{d-\ell}$,
$\ell=1,2,\dots,d-1$.
\item The entries along rows are weakly increasing.
\item The entries along columns are strictly increasing.
\item If at the bottom of column $d-\ell$ we write
$i_d-i_{\ell}-s_\ell+1$, $\ell=1,2,\dots,d-1$ 
(see Figure~8 for the corresponding extended
array in our running example), then rows are still weakly
increasing and columns are still strictly increasing.
\endroster

\proclaim{Corollary 6} The multiplicity 
$M_{\bold j}(\bold i)$ of a point $x\in X_{\bold j}^\circ\subseteq 
X_{\bold i}$ is equal to the number of arrays of positive integers
satisfying {\rm (1)--(4)} above.
\endproclaim

\midinsert
\vskip8pt
\vbox{
$$
\smatrix \format\sa\c\s\c\s\c\s\c\s\c\s\c\s\c\s\c\se\\
\hlinefor{17}\\
&1&&1&&1&&1&&4&&6&&6&&9&\\
\hlinefor{17}\\
&\boldkey 2&&2&&2&&3&&5&& \boldkey1\boldkey0&&\boldkey1\boldkey0 &&11&\\
\hlinefor{17}\\
\omit& && \boldkey3&&4&& \boldkey6&&6&& &\omit& &&\boldkey1\boldkey2&\\
\omit& \omit&\hlinefor9& \omit&\omit& \omit&\hlinefor3\\
\omit& &\omit& &&\boldkey5&& &&\boldkey7&& \\
\omit& \omit&\omit& \omit&\hlinefor3& \omit &\hlinefor3
\endsmatrix
$$
\centerline{\eightpoint Figure 8}
}
\vskip8pt
\endinsert

\enddocument